\documentclass[11pt,reqno]{article}
\usepackage{amsmath,amsthm,amssymb,amscd,epsfig}

\newtheorem{theo}{Theorem}
\newtheorem{lem}[theo]{Lemma}
\newtheorem{prop}[theo]{Proposition}
\newtheorem{coro}[theo]{Corollary}
\newtheorem{exam}[theo]{Example}

\newtheorem*{rem}{Remark}

\def\R{\mathbb R}
\def\C{\mathbb C}
\def\D{\mathbb D}

\def\H{\mathcal H}
\def\M{\mathcal M}

\def\diam{\text{diam}}

\def\Lip{\operatorname{Lip}}

\def\Im{\operatorname{Im}}

\def\supp{\operatorname{supp}}

\title{Beltrami equations with coefficient in the Sobolev space $W^{1,p}$}
\author{A. Clop, D. Faraco, J. Mateu, J. Orobitg, X. Zhong
\thanks{{Clop, Mateu, Orobitg were supported by projects MTM2004-00519, HF2004-0208, 2005-SGR-00774. Faraco was supported by the project Simumat and by Ministerio de Educaci\'on y Ciencia, projects MMM2005-07652-C02-01. Zhong was partially supported by the Academy of Finland, project 207288.
 }\newline
\newline Keywords   Quasiconformal, Hausdorff measure, Removability}}
\date{}

\begin{document}

\maketitle

\begin{abstract}
We study the removable singularities for solutions to the Beltrami
equation $\overline\partial f=\mu\, \partial f$,  assuming that
the coefficient $\mu$ lies on some Sobolev space $W^{1,p}$, $p\leq
2$. Our results are based on an extended version of the well known
Weyl's lemma, asserting that distributional solutions are actually
true solutions. Our main result is that quasiconformal mappings
with compactly supported Beltrami coefficient $\mu\in W^{1,2}$
preserve compact sets of $\sigma$-finite length and vanishing
analytic capacity, even though they need not be bilipschitz.
\end{abstract}

\section{Introduction}

A homeomorphism between planar domains $\phi:\Omega\to\Omega'$ is
called $\mu$-quasiconformal if it is of class
$W^{1,2}_{loc}(\Omega)$ and satisfies the Beltrami equation,
\begin{equation} \label{Beltrami}
\overline\partial \phi(z)=\mu(z)\,\partial\phi(z)
\end{equation}
for almost every $z\in\Omega$. Here $\mu$ is the Beltrami coefficient, that is, a
measurable bounded function with $\|\mu\|_\infty<1$. More
generally, any $W^{1,2}_{loc}(\Omega)$ solution is called
$\mu$-quasiregular. When $\mu=0$, we recover conformal mappings
and analytic functions, respectively.\\
When $\|\mu\|_\infty\leq\frac{K-1}{K+1}<1$ for some $K\geq 1$,
then clearly $\mu$-quasiregular mappings are $K$-quasiregular
 \cite{LV}. As holomorphic mappings are linked to harmonic functions, so
 are quasiregular mappings  with elliptic equations.
 In particular, it is well known that if
 $f=u+iv$ solves (\ref{Beltrami}), $u$ is a solution to

\begin{equation}
\textrm{div}(\sigma \nabla u)=0
\end{equation}
where for almost every $z$, $ (\sigma_{ij}(z)) \in \mathbb{S}(2)$,
the space of symmetric matrices with $\det(\sigma)=1$. Moreover,
$\mu$ and $\sigma$ are related by the formula $\mu=
\frac{\sigma_{11}-\sigma_{22}+2i
\sigma_{12}}{\sigma_{11}+\sigma_{22}+2}$. A similar equation holds
for $v$. Thus it follows from Morrey's work \cite{Mor} that
$\mu$-quasiregular functions are H\"older  continuous with
exponent $\frac{1}{K}$, see also \cite{FS}. On the other hand, if the coefficient $\mu$ is
more regular, say H\"older continuous, then classical Schauder estimates tell us that the derivatives of any $\mu$-quasiregular mapping $f$ are H\"older
continuous and, in particular, $f$ is Lipschitz continuous.\\
\\
In this paper, we study properties of $\mu$-quasiregular mappings (and hence to the solutions to the related elliptic equations), when the regularity of  Beltrami coefficient is measured in the category of Sobolev spaces $W^{1,p}$.  Our study is focussed on the removable singularities and on the distortion of Hausdorff measures and capacities under solutions of such $PDE$.\\
\\
We say that a compact set $E$ is {\it{removable for bounded $\mu$-quasiregular mappings}} if for every open set $\Omega$, every bounded function $f$, $\mu$-quasiregular on $\Omega\setminus E$, admits an extension $\mu$-quasiregular in all of $\Omega$. When $\mu=0$, this is the classical Painlev\'e problem. It is  also a natural question to replace bounded functions by others, such as $BMO$ (bounded mean oscillation),  $VMO$ (vanishing mean oscillation) or $\Lip_\alpha(\Omega)$ (H\"older continuous with exponent $\alpha$). We want to give geometric characterizations of these sets. Of special interest is the case $\mu \in W^{1,2}$, which is at the borderline.\\
\\
By means of the Stoilow factorization (see for instance \cite{LV}), one easily sees that the way $\mu$-quasiconformal mappings distort sets is very related to removability problems. We have the precise bounds for distortion of Hausdorff dimension from Astala \cite{A}, which apply to any $K$-quasiconformal mapping,
$$\dim(\phi(E))\leq\frac{2K\,\dim(E)}{2+(K-1)\,\dim(E)}.$$
In the particular case $\dim(E)=\frac{2}{K+1}$ we also have absolute continuity of measures \cite{ACMOU}, that is,
$$\H^\frac{2}{K+1}(E)=0\hspace{1cm}\Longrightarrow\hspace{1cm}\H^1(\phi(E))=0.$$
If the coefficient is more regular one improves these estimates as well. For instance, if the Beltrami coefficient $\mu$ lies in $VMO$, then every $\mu$-quasiconformal mapping $\phi$ has distributional derivatives in $L^p_{loc}$ for every $p\in(1,\infty)$ (see for instance \cite{AIS} or \cite{GMV}). Thus, $\phi\in\Lip_\alpha$ for every $\alpha\in(0,1)$, and as a consequence, 
$$\dim(\phi(E))\leq\dim(E).$$
Moreover, actually for such $\mu$ one has $\dim(\phi(E))=\dim(E)$.\\
\\
However if we further know that $\mu \in W^{1,2}$ we obtain more precise information. An important reason is the following: We first recall that for $\mu=0$ we have the well known Weyl's Lemma, which asserts that if $T$ is any (Schwartz) distribution such that
$$\langle\overline\partial T, \varphi\rangle=0$$
for each test function $\varphi\in{\cal D}$ (by $\cal D$ we mean the algebra of compactly supported ${\cal C}^\infty$ functions), then $T$ agrees with a holomorphic function. In other words, {\it{distributional}} solutions to Cauchy-Riemann equation are actually {\it{strong}} solutions. When trying to extend this kind of result to the Beltrami equation, one first must define the distribution $(\overline\partial-\mu\,\partial)T=\overline\partial T-\mu\,\partial T$. It need not to make sense, because bounded functions in general do not multiply distributions nicely. However, if the multiplier is asked to exhibit some regularity, and the distribution $T$ is an integrable function, then something may be done. Namely, one can write
$$\langle(\overline\partial-\mu\,\partial)f,\varphi\rangle=-\langle f,\overline\partial\varphi\rangle+\langle f,\partial\mu\,\varphi\rangle+\langle f,\mu\,\partial\varphi\rangle$$
whenever each term makes sense. For instance, this is the case if $\mu\in W^{1,p}_{loc}$ and $f\in L^q_{loc}$, $\frac{1}{p}+\frac{1}{q}=1$. Hence we can call $\overline\partial f-\mu\,\partial f$ the {\it{Beltrami distributional derivative}} of $f$, and we can say that a function $f\in L^q_{loc}$ is {\it{distributionally $\mu$-quasiregular}\index{quasiregular!distributionally quasiregular}} precisely when $(\overline\partial-\mu\,\partial)f=0$ as a distribution. Of course, a priori such functions $f$ could be not quasiregular, since it is not clear if the distributional equation actually implies $f\in W^{1,2}_{loc}$. Thus it is natural to ask when this happens. 

\begin{theo}\label{maintheo6}
Let $f\in L^p_{loc}(\Omega)$ for some $p>2$, and let $\mu\in W^{1,2}$ be a compactly supported Beltrami coefficient. Assume that
$$\langle(\overline\partial-\mu\,\partial)f,\varphi\rangle=0$$
for any $\varphi\in{\cal D}(\Omega)$. Then, $f$ is $\mu$-quasiregular. In particular, $f\in W^{1,2}_{loc}(\Omega)$.
\end{theo}

We must point out that this selfimprovement of regularity is even stronger, because of the factorization theorem for $\mu$-quasiregular mappings, as well as the regularity of homeomorphic solutions when the Beltrami coefficient is {\it{nice}}. Namely, when $\mu\in W^{1,2}$ is compactly supported, it can be shown that any $\mu$-quasiconformal mapping is actually in $W^{2,q}_{loc}$ whenever $q<2$. Hence, every $L^{2+\varepsilon}_{loc}$ distributional solution to the corresponding Beltrami equation is actually a $W^{2,q}_{loc}$ solution, for every $q<2$. Further, we can show that the above Weyl's Lemma holds, as well, when $\mu\in W^{1,p}$ for $p\in(\frac{2K}{K+1},2)$.\\
\\
One may use this selfimprovement to give removability results and study distortion problems for $\mu$-quasiconformal mappings. The conclusions we obtain encourage us to believe that Beltrami equation with $W^{1,2}$ Beltrami coefficient is not so far from the classical planar Cauchy-Riemann equation. For instance, we shall show that for any $0<\alpha<1$, any set $E$ with $\H^{1+\alpha}(E)=0$ is removable for $\Lip_\alpha$ $\mu$-quasiregular mappings, precisely as it is when $\mu=0$ \cite{Do}. Nevertheless, this $\Lip_\alpha$ removability problem does not imply in general any result on $\mu$-quasiconformal distortion of Hausdorff measures, since there are examples of $\mu\in W^{1,2}$ for which the space $\Lip_\alpha$ is not $\mu$-quasiconformally invariant. Therefore, to get results in terms of distortion we study the removability problem with the $BMO$ norm. Then we get that $E$ is removable for $BMO$ $\mu$-quasiregular mappings, if and only if $\H^1(E)=0$. More precisely, this is what happens for $\mu=0$ \cite{K}. Moreover, $E$ is removable for $VMO$ $\mu$-quasiregular mappings if and only if $\H^1(E)$ is $\sigma$-finite, again as in the analytic case \cite{V}. In distortion terms, this reads as $\H^1(E)=0$ if and only if $\H^1(\phi(E))=0$, and $\H^1(E)$ is $\sigma$-finite if and only if $\H^1(\phi(E))$ is.\\
\\
$\mu$-quasiconformal distortion of analytic capacity is somewhat deeper, since the rectifiable structure of sets plays an important role there. We show in Lemma \ref{muqcrectifiabledistortion} that if $\mu\in W^{1,2}$ is compactly supported, then $\phi$ maps rectifiable sets to rectifiable sets. As a consequence, purely unrectifiable sets are mapped to purely unrectifiable sets. Therefore, we get from \cite{D} our following main result.

\begin{theo}\label{maintheo7}
Let $\mu\in W^{1,2}$ be a compactly supported Beltrami coefficient, and let $\phi$ be $\mu$-quasiconformal. If $E$ has $\sigma$-finite length, then
$$\gamma(E)=0\hspace{1cm}\Longleftrightarrow\hspace{1cm}\gamma(\phi(E))=0.$$
\end{theo}

Let's mention that in \cite{T2} Tolsa proved that an homeomorphism $\phi$ preserves the analytic capacity of sets if and only if $\phi$ is a bilipschitz map. On the other hand, the radial stretching $g(z)= z|z|^{\frac{1}{K}-1}$ is not bilipschitz but clearly it preserves sets of zero analytic capacity. Theorem \ref{maintheo7} asserts that $\mu$-quasiconformal mappings , $\mu\in W^{1,2}$, also preserve sets of zero analytic capacity having also $\sigma$-finite length.\\
\\
As a natural question, one may ask wether these distortion results apply also for compactly supported Beltrami coefficients $\mu\in W^{1,p}$ when $p\in(\frac{2K}{K+1},2)$. In this case, we study the same removability problems and we obtain analogous results. For instance, if $\H^{1+\alpha}(E)=0$, then $E$ is removable for $\Lip_\alpha$ $\mu$-quasiregular mappings, as well as for the analytic case \cite{OF}. Again, this does not translate to the distortion problem for Hausdorff measures, since $\Lip_\alpha$ is not quasiconformally invariant. However, this has some interesting consequences in terms of distortion of Hausdorff dimension. Namely, it follows that
\begin{equation}\label{onetoone}
\dim(E)\leq 1\hspace*{1cm}\Longrightarrow\hspace*{1cm}\dim(\phi(E))\leq 1
\end{equation}
Moreover, when letting $\alpha=0$ we get the corresponding $BMO$ and $VMO$ removability problems. Due to our Weyl type Lemma, we show that in this weaker situation $\mu\in W^{1,p}$, $\frac{2K}{K+1}<p<2$, we actually have absolute continuity of measures, i.e.
\begin{equation}\label{zerotozerointro}
\H^1(E)=0\hspace{1cm}\Longrightarrow\hspace{1cm}\H^1(\phi(E))=0
\end{equation}
for any $\mu$-quasiconformal mappig $\phi$. This improves the absolute continuity results in \cite{ACMOU}.\\
\\
We do not know if implication (\ref{zerotozerointro}) is an equivalence. Indeed, if $\frac{2K}{K+1}<p<2$ and $\mu\in W^{1,p}$ then the Beltrami coefficient $\nu$ of inverse mapping $\phi^{-1}$ need not belong to the same Sobolev space $W^{1,p}$ (this is true for $p=2$). However, if $p$ ranges the smaller interval $(\frac{2K^2}{K^2+1},2)$ then a calculation shows that $\nu\in W^{1,r}$ for some $r>\frac{2K}{K+1}$. As a consequence, we obtain the following result.

\begin{theo}\label{maintheo8}
Let $\frac{2K^2}{K^2+1}<p<2$. Let $\mu\in W^{1,p}$ be a compactly supported Beltrami coefficient, and let $\phi$ be $\mu$-quasiconformal. Then,
$$\gamma(E)=0\hspace{1cm}\Longleftrightarrow\hspace{1cm}\gamma(\phi(E))=0,$$
for any compact set $E$ with $\sigma$-finite length.
\end{theo}

In particular, the above result holds whenever our Beltrami coefficient $\mu$ lives in $W^{1,1+\varepsilon}$ and $\|\mu\|_\infty\lesssim\varepsilon$, $\varepsilon>0$.\\
\\
This paper is structured as follows. In Section 2 we study the regularity of $\mu$-quasiregular mappings. In Section 3 we prove Theorem \ref{maintheo6}. In Section 4, we study the $BMO$ and $VMO$ removability problems for $\mu\in W^{1,2}$, and deduce
distortion theorems for $\H^1$. In Section 5 we study $\mu$-quasiconformal distortion of rectifiable sets, and prove Theorem \ref{maintheo7}. In Section 6 we study Beltrami equations with coefficient in $W^{1,p}$, $p\in (\frac{2K}{K+1}, 2)$, and prove Theorem \ref{maintheo8}.

\section{Regularity of $\mu$-quasiconformal mappings}

It is well known (see for instance \cite{BI}) that any $K$-quasiregular mapping belong to better Sobolev spaces than the usual $W^{1,2}_{loc}$ appearing in its definition. More precisely,
$Df\in L^{\frac{2K}{K-1},\infty}$ \cite{A}, and this is sharp. However, if we look not only at
$K$ but also at the regularity of the Beltrami coefficient, something better may be said. This
situation is given when the Beltrami coefficients are in $\Lip_\alpha$. In this case, every homeomorphic solution
(and hence the corresponding $\mu$-quasiregular mappings) have first order derivatives
also in $\Lip_\alpha$. In particular, $\phi$ is locally bilipschitz. The limiting situation in terms of
continuity is obtained when assuming $\mu\in VMO$. In this case, as mentioned before, every $\mu$-quasiconformal mapping has derivatives in $L^p_{loc}(\C)$ for every $p\in(1,\infty)$. Let us discuss the situation in terms of the Sobolev regularity of $\mu$. If $\mu\in W^{1,p}$, $p>2$, then $D\phi\in\Lip_{1-\frac{2}{p}}$, as shows \cite{Ve}. Actually, it comes from \cite{Al} that $\phi\in W^{2,p}_{loc}$. 
In the next lemma we study what happens for an arbitrary $1<p<\infty$.

\begin{prop}\label{regularity}
Let $\mu\in W^{1,p}$ be a compactly supported Beltrami coefficient, and assume
that $\|\mu\|_\infty\leq\frac{K-1}{K+1}$. Let $\phi$ be $\mu$-quasiconformal.
\begin{enumerate}
\item[(a)] If $p>2$, then $\phi\in W^{2,p}_{loc}(\C)$.
\item[(b)] If $p=2$, then $\phi\in W^{2,q}_{loc}(\C)$ for every $q<2$.
\item[(c)] If $\frac{2K}{K+1}<p<2$, then $\phi\in W^{2,q}_{loc}(\C)$ for every $q<q_0$,
where $\frac{1}{q_0}=\frac{1}{p}+\frac{K-1}{2K}$.
\end{enumerate}
\end{prop}
\begin{proof}
There is no restriction if we suppose that $\mu$ has compact support included in $\D$.
Assume first that $p>2$. Then, arguing as in Ahlfors \cite[Lemma 5.3]{Al}, there exists a continuous function $g$ such that $\partial\phi=e^g$. This function $g$ is a solution to $\overline\partial g=\mu\,\partial g + \partial\mu$, which may be constructed as
$$g=\frac{1}{z}\ast(I-\mu B)^{-1}(\partial\mu)$$
where $B$ denotes the Beurling transform. Clearly, $g\in W^{1,p}_{loc}(\C)$ and, since it is continuous, also $e^g$ is contiuous. Moreover, $\partial(e^g)=e^g\,\partial g$ and the same happens with $\overline\partial$. Then, $\partial\phi\in W^{1,p}_{loc}(\C)$ and the result follows.\\
\\
Let now $p\le 2$. Let $\psi\in{\cal C}^\infty(\C)$, $0\leq\psi\leq 1$, $\int\psi=1$, supported on $\D$, and let $\psi_n(z)=n^2\,\psi(nz)$. Define $$\mu_n(z)=\mu\ast\psi_n(z)=\int n^2\,\psi(nw)\,\mu(z-w)\,dA(w)$$
Then $\mu_n$ is of class ${\cal C}^\infty$, has compact support inside of $2\D$, $\|\mu_n\|_\infty\leq\|\mu\|_\infty$ and $\mu_n\to \mu$ in $W^{1,p}(\C)$ as $n\to\infty$. As in \cite{Al}, the corresponding principal solutions $\phi_n$ and $\phi$ can be written as $\phi(z)=z+\tilde{C}h(z)$ and $\phi_n(z)=z+\tilde{C}h_n(z)$, where $h, h_n$ are respectively defined by $h=\mu Bh + \mu$ and $h_n=\mu_n Bh_n+ \mu_n$. We then get $\phi_n\to\phi$ as $n\to\infty$ with convergence in $W^{1,r}$ for every $r<\frac{2K}{K-1}$. Now observe that $\phi_n$ is a ${\cal C}^\infty$ diffeomorphism and conformal outside of $2\D$. This allows us to take derivatives in the equation $\overline\partial\phi_n=\mu_n\,\partial\phi_n$. We get
$$\overline\partial\partial\phi_n-\mu_n\,\partial\partial\phi_n=\partial\mu_n\,\partial\phi_n.$$
This may be written as
$$(\overline\partial-\mu_n\,\partial)(\log\partial\phi_n)=\partial\mu_n$$
or equivalently
\begin{equation}\label{logjacobian}
(I-\mu_n B)(\overline\partial\log(\partial\phi_n))=\partial\mu_n
\end{equation}
so that
\begin{equation}\label{laplacian}
\overline\partial\partial\phi_n=\partial\phi_n\,(I-\mu_nB)^{-1}(\partial\mu_n)
\end{equation}
Fix $\frac{2K}{K+1}<p<2$. In this case \cite{AIS}, the norm of $\|(I-\mu_n B)^{-1}\|_{L^p(\C) \to L^p(\C)}$ depends only on $K$ and $p$. Now recall that $\partial\mu_n\to\partial\mu$ in $L^p(\C)$ and $\partial \phi_n \to \partial \phi$ in  $L^r(\C)$ for $r<\frac{2K}{K-1}$. Then  if $q<q_0$, $\frac{1}{q_0}=\frac{1}{p}+\frac{K-1}{2K}$, the right hand side in (\ref{laplacian}) converges to $(I-\mu B)^{-1}(\partial\mu)\,\partial\phi$ in $L^q(\C)$. Hence, the sequence $(\overline\partial\partial\phi_n)_n$ is uniformly bounded in $L^q(\C)$. Taking a subsequence, we get that $\phi_n$ converges in $W^{2,q}_{loc}(\C)$, and obviously the limit is $\phi$, so that $\phi\in W^{2,q}_{loc}(\C)$.\\
\\
Assume finally that $p=2$. Repeating the argument above, we get $\phi\in W^{2,q}$ for every $q<\frac{2K}{2K-1}<2$, which is weaker than the desired result. To improve it, we first show that $\phi_n\to\phi$ in $W^{1,r}_{loc}(\C)$ for every $r\in(1,\infty)$. To do that, notice that both $I-\mu_n B$ and $I-\mu B$ are invertible operators in $L^r(\C)$ for all $r\in(1,\infty)$, since both $\mu_n, \mu\in VMO$ (see for instance \cite{I}). Further, from the Sobolev Embedding Theorem, $\mu_n\to \mu$ in $L^r(\C)$. Thus,
$$\lim_{n\to\infty}\|(I-\mu_n B)-(I-\mu B)\|_{L^r\to L^r}=\lim_{n\to\infty}\|(\mu_n-\mu)B\|_{L^r\to L^r}=0$$
for any $r\in(1,\infty)$. Now recall that the set of bounded operators $L^r(\C)\to L^r(\C)$ defines a complex Banach algebra, in which the invertible operators are an open set and, moreover, the inversion is continuous. As a consequence,
$$\lim_{n\to\infty}\|(I-\mu_n B)^{-1}\|_{L^r\to L^r} = \|(I-\mu B)^{-1}\|_{L^r\to L^r}$$
for each $r\in(1,\infty)$. This implies that $h_n\to h$ in $L^r(\C)$ so that $\phi_n\to\phi$ in $W^{1,r}_{loc}(\C)$. Going back to (\ref{laplacian}), the right hand side converges to $\partial\phi\,(I-\mu B)^{-1}(\partial\mu)$ in the norm of $L^q(\C)$, provided that $q<2$, and now the result follows.
\end{proof}

If $p>2$, $D^2\phi$ cannot have better integrability than $D\mu$, since in that case $J(\cdot,\phi)$ is a continuous function, bounded from above and from below. If $p=2$, the sharpness of the above proposition may be stated as a consequence of the following example \cite[p.142]{Va}.

\begin{exam}\label{W12example}
The function $\phi(z)=z\,(1-\log|z|)$ is $\mu$-quasiconformal in a neighbourhood of the origin, with Beltrami coefficient
\begin{equation}
\mu(z)=\frac{z}{\overline{z}}\,\frac{1}{2\log|z|-1}
\end{equation}
In particular, we have $\mu\in W^{1,2}$ in a neighbourhood of the origin. Thus, we have $\phi\in W^{2,q}_{loc}$ whenever $q<2$. However,
$$|D^2\phi(z)|\simeq\frac{1}{|z|}$$
so that $\phi\notin W^{2,2}_{loc}$.
\end{exam}

Finally, the radial stretching $f(z)=z|z|^{\frac{1}{K}-1}$ has Beltrami coefficient in $W^{1,p}$ for every $p<2$ and, however, $D^2f$ lives in no better space than $L^{\frac{2K}{2K-1},\infty}$.\\
\\
In order to study distortion results, we need information about the integrability of the inverse of a $\mu$-quasiconformal mapping. This can be done by determining the Sobolev regularity of the corresponding Beltrami coefficient to $\phi^{-1}$.

\begin{prop}\label{muinverse}
Let $\mu\in W^{1,2}$ be a compactly supported Beltrami coefficient, and let $\phi$ be $\mu$-quasiconformal. Then, $\phi^{-1}$ has Beltrami coefficient
$$\nu(z)=-\mu(\phi^{-1}(z))\,\frac{\partial\phi}{\overline{\partial\phi}}(\phi^{-1}(z))$$
In particular, $\nu\in W^{1,2}$.
\end{prop}
\begin{proof}
An easy computation shows that
$$
\nu(z)=\frac{\overline\partial\phi^{-1}(z)}{\partial\phi^{-1}(z)}=-\left(\mu\,\frac{\partial\phi}{\overline{\partial\phi}}\right)(\phi^{-1}(z))
$$
For compactly supported $\mu\in W^{1,2}$, it follows from equation (\ref{logjacobian}) that the normalized solution $\phi$ is such that $\log\partial\phi\in W^{1,2}$. Hence,
$$\partial\phi=e^\lambda$$
for a function $\lambda\in W^{1,2}(\C)$ (in fact, $\lambda=\log\partial\phi$).  
Thus, in terms of $\lambda$, we get
$$\nu\circ\phi=-\mu\,e^{2i\,\Im(\lambda)}$$
where $\Im(\lambda)$ is the imaginary part of the function $\lambda$. Hence,
$$D(\nu\circ\phi)=-D\mu\,e^{2i\,\text{Im}(\lambda)}-\mu\,2i\,e^{2i\,\text{Im}(\lambda)}\,D(\text{Im}(\lambda))$$
so that
$$|D(\nu\circ\phi)|\leq |D\mu|+2\,|\mu|\,|D(\text{Im}(\lambda))|$$
In particular, $\nu\circ\phi$ has derivatives in $L^2(\C)$. 
Now, from the identity
$$\int|D\nu(z)|^2dA(z)=\int|D\nu(\phi(w))|^2J(w,\phi)\,dA(w)\leq\int|D(\nu\circ\phi)(w)|^2dA(w)$$
the result follows. 
\end{proof}

\begin{rem} \label{inversecoefficient}
As shown above, if $\mu$ belongs to $W^{1,p}$ for some $p \ge 2$, then
the same can be said for $\nu$. If $\mu$ is only in $W^{1,p}$,
$\frac{2K}{K+1}<p<2$, the situation is different. More precisely,
an argument as above shows that $\nu\in W^{1,r}$
for every $r$ such that
$$r<\frac{2p}{2K-(K-1)p}$$
In particular, for $p>\frac{2K}{K+1}$ we always have $\nu\in
W^{1,1}$, but $\nu$ does not fall, in general, in the same Sobolev
space $W^{1,p}$ than $\mu$. However, for $p>\frac{2K^2}{K^2+1}$,
we always have $\nu\in W^{1,r}$ for some $r>\frac{2K}{K+1}$.
\end{rem}

The above regularity results can be applied to study distortion properties of $\mu$-quasiconformal mappings. For instance, if $\mu$ is a compactly supported $W^{1,2}$ Beltrami coefficient, then both $\phi$ and $\phi^{-1}$ are $W^{2,q}_{loc}$ functions, for every $q<2$. Therefore, $\phi,\phi^{-1}\in\Lip_\alpha$ for every $\alpha\in(0,1)$ (notice that this is true under the more general assumption $\mu\in VMO$). Thus,
\begin{equation}\label{fixdimension}
\dim(\phi(E))=\dim(E)
\end{equation}
On the other hand, we may ask if this identity can be translated to Hausdorff measures. As a matter of fact, observe that the mapping in Example \ref{W12example} is not Lipschitz continuous. Thus, is not clear how $\mu$-quasiconformal mappings with $W^{1,2}$ Beltrami coefficient distort Hausdorff measures or other set functions, such as analytic capacity, even preserving Hausdorff dimension. Further, we do not know if for Beltrami coefficients $\mu\in W^{1,p}$, $p<2$, the corresponding $\mu$-quasiconformal mappings satisfy equation (\ref{fixdimension}) or not. Questions related with this will be treated in Sections \ref{muqcdistortion}, \ref{boundedlymuremovablesets} and \ref{muW1p}.

\section{Distributional Beltrami equation with $\mu\in W^{1,2}$}

A typical feature in the theory of quasiconformal mappings is the
selfimprovement of regularity. Namely, it is known that
weakly $K$-quasiregular mappings in $W^{1,\frac{2K}{K+1}}_{loc}$
are actually $K$-quasiregular \cite{A, PV}. This improvement is stronger for
$K=1$, since in this case we do not need any Sobolev regularity as
a starting point. The classical Weyl's Lemma establishes that if
$f$ is a distribution such that
$$\langle\overline\partial f, \varphi\rangle=0$$
for every testing function $\varphi\in{\cal D}$, then $f$ agrees at almost every point
with a holomorphic function. Our following goal is to deduce an extension to this result
for the Beltrami operator, provided that $\mu\in W^{1,2}$.\\
\\
Assume first we are given a compactly supported Beltrami coefficient $\mu$, such that $\mu\in W^{1,2}$. Let
$f\in L^p_{loc}$ for some $p\in(2,\infty)$. We can define a linear functional
$$\langle\overline\partial f- \mu\,\partial f,\varphi\rangle=-\langle f,(\overline\partial-\partial\mu)\varphi\rangle=-\langle f,\overline\partial\varphi\rangle+\langle f,\partial(\mu\,\varphi)\rangle$$
for each compactly supported $\varphi\in{\cal C}^\infty$. Clearly, $\overline\partial f-\mu\,\partial f$ defines a distribution, which will be called the {\it{Beltrami distributional derivative}} of $f$.\\
\\
We say that a function $f\in L^p_{loc}$ is {\it{distributionally $\mu$-quasiregular}} if its Beltrami distributional derivative vanishes, that is,
$$\langle\overline\partial f-\mu\,\partial f,\varphi\rangle=0$$
for every testing function $\varphi\in{\cal D}$. It turns out that one may take then a bigger class of testing functions $\varphi$.

\begin{lem}\label{previous2}
Let $p>2$, $q=\frac{p}{p-1}$, and let $\mu\in W^{1,2}$ be a compactly supported Beltrami coefficient. Assume that $f\in L^p_{loc}$ satisfies
$$\langle\overline\partial f-\mu\,\partial f, \varphi\rangle=0$$
for every $\varphi\in{\cal D}$. Then, it also holds for $\varphi\in W^{1,q}_0$.
\end{lem}
\begin{proof}
When $\mu\in W^{1,2}$ is compactly supported and $f\in L^p_{loc}$ for some $p>2$, the Beltrami distributional derivative $\overline\partial f-\mu\,\partial f$ acts continuously on $W^{1,q}_0$ functions, since
$$\aligned
|\langle\overline\partial f- \mu\,\partial f,\varphi\rangle|
&\leq|\langle f,\overline\partial\varphi\rangle|+|\langle f,\partial(\mu\,\varphi)\rangle|\\
&\leq\int |f|\,|\overline\partial\varphi|+\int |f|\,|\partial\mu|\,|\varphi|+\int |f|\,|\mu|\,|\partial\varphi|\\
&\leq \|f\|_p\,\|\overline\partial\varphi\|_q+ \|f\|_p\,\|\partial\mu\|_2\,\|\varphi\|_\frac{2q}{2-q}+ \|f\|_p\,\|\mu\|_\infty\,\|\partial\varphi\|_q\\
\endaligned$$
Hence, if $\overline\partial f-\mu\,\partial f$ vanishes when acting on ${\cal D}$, it will also vanish on $W^{1,q}_0$.
\end{proof}

\begin{theo}\label{weylslemmaforbeltrami}
Let $f\in L^p_{loc}$ for some $p>2$. Let $\mu\in W^{1,2}$ be a compactly supported Beltrami coefficient. Assume that
$$\langle\overline\partial f-\mu\,\partial f, \psi\rangle=0$$
for each $\psi\in{\cal D}$. Then, $f$ is $\mu$-quasiregular.
\end{theo}
\begin{proof}
Let $\phi$ be any $\mu$-quasiconformal mapping, and define $g=f\circ\phi^{-1}$. Since $\phi\in W^{2,q}_{loc}$ for any $q<2$, then $J(\cdot,\phi)\in L^q_{loc}$ for every $q\in(1,\infty)$ so that $g\in L^{p-\varepsilon}_{loc}$ for every $\varepsilon>0$. Thus, we can define $\overline\partial g$ as a distribution. We have for each $\varphi\in{\cal D}$
$$\aligned
\langle\overline\partial g,\varphi\rangle
&=-\langle g,\overline\partial\varphi\rangle\\
&=-\int g(w)\,\overline\partial\varphi(w)\,dA(w)\\
&=-\int f(z)\,\overline\partial\varphi(\phi(z))\,J(z,\phi)\,dA(z)\\
&=-\int f(z)\,\big(\partial\phi(z)\,\overline\partial(\varphi\circ\phi)(z)-\overline\partial\phi(z)\,\partial(\varphi\circ\phi)(z)\big)\,dA(z).
\endaligned$$
On one hand,
$$\aligned
-\int f(z)\,\overline\partial\phi(z)\,\partial(\varphi\circ\phi)(z)dA(z)
&=\langle \partial f,\overline\partial\phi\cdot\,\varphi\circ\phi\rangle + \int f(z)\,\partial\overline\partial\phi(z)\,\varphi\circ\phi(z)\,dA(z).\\
\endaligned$$
and here everything makes sense. On the other hand,
$$
-\int f(z)\,\partial\phi(z)\,\overline\partial(\varphi\circ\phi)(z)\,dA(z)=\langle\overline\partial f,\partial\phi\cdot\varphi\circ\phi\rangle + \int f(z)\,\overline\partial\partial\phi(z)\,\varphi\circ\phi(z)\,dA(z).
$$
Therefore,
$$
\langle\overline\partial g,\varphi\rangle
=\langle\overline\partial f,\partial\phi\cdot\varphi\circ\phi\rangle -\langle\partial f,\overline\partial\phi\cdot\varphi\circ\phi\rangle.$$
But if $\varphi\in{\cal D}$ then the function $\psi=\partial\phi\cdot\varphi\circ\phi$ belongs to $W^{1,q}_0$ for every $q<2$ and, in particular, for $q=\frac{p}{p-1}$, provided that $p>2$. Hence, also $\mu\psi\in W^{1,q}$. Thus,
$$\aligned
\langle\overline\partial g,\varphi\rangle
&=\langle\overline\partial f,\partial\phi\,\,\varphi\circ\phi\rangle -\langle\partial f,\mu\,\partial\phi\,\,\varphi\circ\phi\rangle\\
&=\langle\overline\partial f,\partial\phi\,\,\varphi\circ\phi\rangle -\langle\mu\,\partial f,\partial\phi\,\,\varphi\circ\phi\rangle\\
&=\langle\overline\partial f-\mu\,\partial f,\partial\phi\,\,\varphi\circ\phi\rangle.
\endaligned
$$
By Lemma \ref{previous2}, the last term vanishes. Hence, $g$ is holomorphic and therefore $f$ is $\mu$-quasiregular.
\end{proof}

\begin{rem}
It should be said that the argument used in the proof does not work for the generalized Beltrami equation $\overline\partial f=\mu\,\partial f+\nu\,\overline{\partial f}$ (with compactly supported $\mu,\nu\in W^{1,2}$ with $\||\mu|+|\nu|\|_\infty<1$), because in this more general setting there is not Stoilow's factorization theorem. For more information about this equation we refer the reader to \cite{Re}.
\end{rem}

From the above theorem, if $f$ is an $L^p_{loc}$ function for some $p>2$ whose Beltrami distributional derivative vanishes, then $f$ may be written as $f=h\circ\phi$ with holomorphic $h$ and $\mu$-quasiregular $\phi$. As a consequence, we get $f\in W^{2,q}_{loc}$ for every $q<2$, so we actually gain not $1$ but $2$ degrees of regularity.

\section{$\mu$-quasiconformal distortion of Hausdorff measures}\label{muqcdistortion}

Let $E$ be a compact set, and let $\mu$ be any compactly supported
$W^{1,2}$ Beltrami coefficient. If $\phi$ is
$\mu$-quasiconformal, then it follows already from the fact that
$\mu \in VMO$ that $\dim(\phi(E))=\dim(E)$. However we do not
know how Hausdorff measures are distorted.  In this section we
answer this question when $\dim(E)=1$, but in an indirect way.
Our arguments go through some removability problems for
$\mu$-quasiregular mappings. For solving these problems, the
Weyl's Lemma for the Beltrami equation (Theorem
\ref{weylslemmaforbeltrami}) plays an important role.\\
\\
Given a compact set $E$ and two real numbers $t\in(0,2)$ and $\delta>0$, we denote 
$$\M^t_\delta(E)=\inf\left\{\sum_j\diam(D_j)^t; E\subset\cup_j D_j, \diam(D_j)\leq\delta\right\}.$$ 
Then, $\M^t(E)=\M^t_\infty(E)$ is the $t$-dimensional Hausdorff content of $E$, and 
$$\H^t(E)=\lim_{\delta\to 0}\M^t_\delta(E)$$
is the $t$-dimensional Hausdorff measure of $E$. Analogously, for any nondecreasing function $h:[0,\infty)\to [0,\infty)$ with $h(0)=0$, we denote
$$\M^h_\delta(E)=\inf\left\{\sum_jh(\diam(D_j))^t; E\subset\cup_j D_j, \diam(D_j)\leq\delta\right\}.$$ 
and $\M^h(E)=\M^h_\infty(E)$. Then, 
$$
\M^t_\ast(E)=\sup\left\{\M^h(E);\lim_{s\to 0}\frac{h(s)}{s^t}=0\right\}
$$
is called the $t$-dimensional lower Hausdorff content of $E$.

\begin{lem}\label{bmomuremovablesets}
Let $E$ be a compact set, and $\mu\in W^{1,2}$ a Beltrami coefficient,
with compact support inside of $\D$. Suppose that $f$ is $\mu$-quasiregular
 on $\C\setminus E$, and $\varphi\in{\cal D}$.
\begin{enumerate}
\item[(a)] If $f\in BMO(\C)$, then
$$|\langle\overline\partial f-\mu\,\partial f,\varphi\rangle|\leq\,C\,\left(1+\|\mu\|_\infty+\|\partial\mu\|_2\right)\,\left(\|\varphi\|_\infty+\|D\varphi\|_\infty\right)\,\|f\|_\ast\,\M^1(E)$$
\item[(b)] If $f\in VMO(\C)$, then
$$|\langle\overline\partial f-\mu\,\partial f,\varphi\rangle|\leq\,C\,\left(1+\|\mu\|_\infty+\|\partial\mu\|_2\right)\,\left(\|\varphi\|_\infty+\|D\varphi\|_\infty\right)\,\|f\|_\ast\,\M^1_\ast(E)$$
\item[(c)] If $f\in\Lip_\alpha(\C)$, then
$$|\langle\overline\partial f-\mu\,\partial f,\varphi\rangle|\leq\,C\,\left(1+\|\mu\|_\infty+\|\partial\mu\|_2\right)\,\left(\|\varphi\|_\infty+\|D\varphi\|_\infty\right)\,\|f\|_\alpha\,\M^{1+\alpha}(E)$$
\end{enumerate}
\end{lem}
\begin{proof}
We consider the function $\delta=\delta(t)$ defined by
$$\delta(t)=\sup_{\diam(D)\leq 2t}\left(\frac{1}{|D|}\int_D|f-f_D|^2\right)^\frac{1}{2}$$
when $0<t<1$, and $\delta(t)=1$ if $t\geq 1$. By construction, for each disk $D\subset\C$ we have
$$\left(\frac{1}{|D|}\int_D|f-f_D|^2\right)^\frac{1}{2}\leq\delta\left(\frac{\diam(D)}{2}\right).$$
Now consider the measure function $h(t)=t\,\delta(t)$. Let
$\{D_j\}_{j=1}^n$ be a covering of $E$ by disks, such
that
$$\sum_jh(\diam(D_j))\leq\M^h(E)+\varepsilon$$
and consider a partition of unity $\psi_j$ subordinated to the covering $D_j$. Each
$\psi_j$ is a ${\cal C}^\infty$ function, compactly supported in
$2D_j$, $|D\psi_j(z)|\leq\frac{C}{\diam(2D_j)}$ and $0\leq\sum_j\psi_j\leq 1$ on $\C$. In particular, $\sum_j\psi_j=1$ on $\cup_jD_j$. Since $f$ is $\mu$ quasiregular on $\C \setminus E$, we have that for every test function $\varphi\in{\cal D}$,
\begin{equation}\label{dossumes}
-\langle\overline\partial f-\mu\,\partial f,\varphi\rangle=\sum_{j=1}^n\langle f-c_j,\overline\partial(\varphi\,\psi_j)\rangle-\sum_{j=1}^n\langle f-c_j,\partial(\mu\,\varphi\,\psi_j)\rangle.
\end{equation}
For the first sum, we have
$$\aligned
\left|\sum_{j=1}^n\langle f-c_j,\overline\partial(\varphi\,\psi_j)\rangle\right|
&\leq\sum_j\int_{2D_j}|f-c_j|\,\left(|\overline\partial\varphi|\,|\psi_j|+|\varphi|\,|\overline\partial\psi_j|\right)\\
&\leq \sum_j\left(\int_{2D_j}|f-c_j|^2\right)^\frac{1}{2}\,
\left(\|\overline\partial\varphi\|_\infty\diam(2D_j)+C\|\varphi\|_\infty\right)\\
&\lesssim\sum_jh(\diam(D_j))\,\left(\|\overline\partial\varphi\|_\infty
\diam(2D_j)+C\|\varphi\|_\infty\right)\endaligned$$
and this sum may be bounded by $\left(\M^h(E)+\varepsilon\right)\left(\|\varphi\|_\infty+\|D\varphi\|_\infty\right)$. The second sum in (\ref{dossumes}) is divided into two terms,
$$
\left|\sum_{j=1}^n\langle f-c_j,\partial(\mu\,\varphi\,\psi_j)\rangle\right|
\leq\sum_j\int_{2D_j}|f-c_j|\,|\partial\mu|\,|\varphi\,\psi_j|+\sum_j\int_{2D_j}|f-c_j|\,|\mu|\,|\partial(\varphi\,\psi_j)|.$$
The second term can be bounded as before,
$$
\sum_j\int_{2D_j}|f-c_j|\,|\mu|\,|\partial(\varphi\,\psi_j)|\lesssim\|\mu\|_\infty\,\left(\M^h(E)+\varepsilon\right)\left(\|\varphi\|_\infty+\|D\varphi\|_\infty\right).
$$
Finally, for the first term, and using that $0\leq\sum_j\psi_j\leq 1$,
$$\aligned
\sum_j\int|f-c_j|\,|\partial\mu|\,|\varphi\,\psi_j|
&\leq\|\varphi\|_\infty\,\sum_j\left(\int|f-c_j|^2\,|\psi_j|\right)^\frac{1}{2}\,\left(\int|\partial\mu|^2\,|\psi_j|\right)^\frac{1}{2}\\
&\leq\|\varphi\|_\infty\,\left(\sum_j|2D_j|\delta(\diam(D_j))^2\right)^\frac{1}{2}\,
\left(\sum_j\int|\partial\mu|^2\,\psi_j\right)^\frac{1}{2}\\
&\lesssim\|\varphi\|_\infty\,\left(\sum_j\diam(D_j)^2\,\delta(\diam(D_j))^2\right)^\frac{1}{2}\,\left(\int_\D|\partial\mu|^2\right)^\frac{1}{2}\\
&\lesssim\|\varphi\|_\infty\,\|\partial\mu\|_2\,\left|\cup_jD_j\right|^\frac{1}{2}.
\endaligned$$
Observe that this term is harmless since the area $|\cup_j D_j|$ is bounded by
$$
|\cup_j D_j|\lesssim\M^h\left(\cup_jD_j\right)\leq\sum_jh\left(\diam(D_j)\right)\leq\M^h(E)+\varepsilon
$$


It just remains to distinguish in terms of the regularity of $f$. If $f\in BMO(\C)$ then the best we can say is that $\delta(t)\lesssim\|f\|_\ast$ for all $t>0$, so that $\M^h(E)\leq\M^1(E)\|f\|_\ast $. Secondly, if $f\in VMO(\C)$ then we also have
$$\lim_{t\to 0}\frac{h(t)}{t}=\lim_{t\to 0}\delta(t)=0$$
and hence $\M^h(E)\leq\M^1_\ast(E)\|f\|_\ast$. Finally, if $f\in\Lip_\alpha$, then $\delta(t)\lesssim\|f\|_\alpha t^\alpha$, and therefore $\M^h(E)\leq\M^{1+\alpha}(E)\|f\|_\alpha$.
\end{proof}

Lemma \ref{bmomuremovablesets} has very interesting consequences, related to $\mu$-quasiconformal distortion. First, we show that our $\mu$-quasiconformal mappings preserve sets of zero length.

\begin{coro}\label{muzerotozero}
Let $E\subset\C$ be a compact set. Let $\mu\in W^{1,2}$ be a compactly supported Beltrami coefficient, and $\phi$ a $\mu$-quasiconformal mapping. Then,
$$\H^1(E)=0\hspace*{1cm}\Longleftrightarrow\hspace*{1cm}\H^1(\phi(E))=0$$
\end{coro}
\begin{proof}
By Proposition \ref{muinverse}, it will suffice to prove
that $\H^1(E)=0$ implies $\H^1(\phi(E))=0$. Assume, thus, that
$\H^1(E)=0$. Let $f\in BMO(\C)$ be holomorphic on
$\C\setminus\phi(E)$. Then $g=f\circ\phi$ belongs also to
$BMO(\C)$. Moreover, $g$ is $\mu$-quasiregular on $\C\setminus E$
so that, by Lemma \ref{bmomuremovablesets},
$\langle\overline\partial g-\mu\,\partial g,\varphi\rangle=0$
whenever $\varphi\in{\cal D}$. As a consequence, by Theorem
\ref{weylslemmaforbeltrami}, $g$ is $\mu$-quasiregular on the
whole of $\C$ and hence $f$ admits an entire extension. This says
that the set $\phi(E)$ is removable for $BMO$ holomorphic
functions. But these sets are characterized \cite{K} by the
condition $\H^1(\phi(E))=0$.
\end{proof}

Another consequence is the complete solution of the removability problem for $BMO$ $\mu$-quasiregular mappings. Recall that a compact set $E$ is said {\it{removable for $BMO$ $\mu$-quasiregular mappings}} if every function $f\in BMO(\C)$, $\mu$-quasiregular on $\C\setminus E$, admits an extension which is $\mu$-quasiregular on $\C$.

\begin{coro}
Let $E\subset\C$ be compact. Let $\mu\in W^{1,2}$ be a compactly supported Beltrami coefficient. Then, $E$ is removable for $BMO$ $\mu$-quasiregular mappings if and only if $\H^1(E)=0$.
\end{coro}
\begin{proof}
Assume first that $\H^1(E)=0$, and let $f\in BMO(\C)$ be
$\mu$-quasiregular on $\C\setminus E$. Then, by Lemma
\ref{bmomuremovablesets}, we have $\langle\overline\partial
f-\mu\,\partial f,\varphi\rangle=0$ for every $\varphi\in{\cal
D}$. Now by Theorem \ref{weylslemmaforbeltrami} we deduce that $f$
is $\mu$-quasiregular. Consequently, $E$ is removable. Conversely,
if $\H^1(E)>0$, then by Corollary \ref{muzerotozero},
$\H^1(\phi(E))>0$, so that $\phi(E)$ is not removable for $BMO$
analytic functions \cite{K}. Hence, there exists a function $h$
belonging to $BMO(\C)$, holomorphic on $\C\setminus\phi(E)$, non
entire. But therefore $h\circ\phi$ belongs to $BMO(\C)$, is
$\mu$-quasiregular on $\C\setminus E$, and does not admit any
$\mu$-quasiregular extension on $\C$. Consequently, $E$ is not
removable for $BMO$ $\mu$-quasiregular mappings.
\end{proof}

A second family of consequences of Lemma \ref{bmomuremovablesets}
comes from the study of the $VMO$ case. First, we prove that
$\mu$-quasiconformal mappings preserve  compact sets with
$\sigma$-finite length.

\begin{coro}\label{mufinitetofinite}
Let $\mu\in W^{1,2}$ be a compactly supported Beltrami coefficient, and $\phi$ any $\mu$-quasiconformal mapping. For every compact set $E$,
$$\H^1(E)\text{ is $\sigma$-finite}\hspace{1cm}\Longleftrightarrow\hspace{1cm}\H^1(\phi(E))\text{ is $\sigma$-finite}$$
\end{coro}
\begin{proof}
Again, we only have to show that $\M^1_\ast(E)=0$ implies $\M^1_\ast(\phi(E))=0$. Assume, thus, that $\M^1_\ast(E)=0$, and let $f\in VMO(\C)$ be analytic on $\C\setminus\phi(E)$. If we prove that $f$ extends holomorphically on $\C$, then $\phi(E)$ must have $\sigma$-finite length, and we will be done. To do that, we first observe that $g=f\circ\phi$ also belongs to $VMO(\C)$. Further, $g$ is $\mu$-quasiregular on $\C\setminus E$, and since $\M^1_\ast(E)=0$, by Lemma \ref{bmomuremovablesets} we get that $\overline\partial g-\mu\,\partial g=0$ on ${\cal D}'$. Consequently, from Theorem \ref{weylslemmaforbeltrami}, $g$ is $\mu$-quasiregular on the whole of $\C$ and hence $f$ extends holomorphically on $\C$.
\end{proof}

As in the $BMO$ setting, the removability problem for $VMO$ $\mu$-quasiregular functions also gets solved. A compact
 set $E$ is said to be {\it{removable for $VMO$ $\mu$-quasiregular mappings}} if every function
 $f\in VMO(\C)$ $\mu$-quasiregular on $\C\setminus E$ admits an extension which is $\mu$-quasiregular on $\C$.

\begin{coro}
Let $E\subset\C$ be compact. Let $\mu\in W^{1,2}$ be a compactly supported Beltrami coefficient.
Then $E$ is removable for $VMO$ $\mu$-quasiregular mappings if and only if $\H^1(E)$ is $\sigma$-finite.
\end{coro}
\begin{proof}
If $\H^1(E)$ is $\sigma$-finite, then $\M^1_\ast(E)=0$, so that from Lemma \ref{bmomuremovablesets} every
 function $f\in VMO(\C)$ $\mu$-quasiregular on $\C\setminus E$ satisfies
  $\overline\partial f=\mu\,\partial f$ on ${\cal D}'$. By Theorem \ref{weylslemmaforbeltrami}, $f$
  extends $\mu$-quasiregularly and thus $E$ is removable.\\
If $\H^1(E)$ is not $\sigma$-finite, we have just seen that $\H^1(\phi(E))$ must not be $\sigma$-finite.
Thus, it comes from Verdera's work \cite{V} that there exists a function $h\in VMO(\C)$, analytic
on $\C\setminus\phi(E)$, non entire. But therefore $h\circ\phi$ belongs to $VMO$, is $\mu$-quasiregular on
$\C\setminus E$, and does not extend $\mu$-quasiregularly on $\C$.
\end{proof}

The class $\Lip_\alpha$ has, in comparison with $BMO$ or $VMO$, the disadvantage
of being not quasiconformally invariant. This means that we cannot read any removability result
for $\Lip_\alpha$ in terms of distortion of Hausdorff measures, and therefore for $\H^{1+\alpha}$ we
cannot obtain results as precise as Lemmas \ref{muzerotozero} or \ref{mufinitetofinite}. 
Hence question remains unsolved.
%
%
However, Theorem \ref{weylslemmaforbeltrami} can be used to study the $\Lip_\alpha$ removability problem. Recall that a compact set $E$ is {\it{removable for $\Lip_\alpha$ $\mu$-quasiregular mappings}} if every function $f\in\Lip_\alpha(\C)$, $\mu$-quasiregular on $\C\setminus E$, has a $\mu$-quasiregular extension on $\C$.

\begin{coro}
Let $E$ be compact, and assume that $\H^{1+\alpha}(E)=0$. Then, $E$ is removable for $\Lip_\alpha$ $\mu$-quasiregular mappings.
\end{coro}
\begin{proof}
As before, if $f\in\Lip_\alpha$ is $\mu$-quasiregular outside of $E$, then Lemma \ref{bmomuremovablesets} tells us that its Beltrami distributional derivative vanishes. By Theorem \ref{weylslemmaforbeltrami}, we get that $f$ is $\mu$-quasiregular.
\end{proof}

The above result is sharp, in the sense that if $\H^{1+\alpha}(E)>0$ then there is a compactly supported Beltrami coefficient $\mu\in W^{1,2}$ such that $E$ is not removable for $\Lip_\alpha$ $\mu$-quasiregular mappings (take simply $\mu=0$, \cite{OF}). 
%

\section{$\mu$-quasiconformal distortion of analytic capacity}\label{boundedlymuremovablesets}

If $\mu\in W^{1,2}(\C)$ is a Beltrami coefficient, compactly supported on $\D$, and $E\subset\D$ is compact, we say that $E$ is {\it{removable for bounded $\mu$-quasiregular functions}}, if and only if any bounded function $f$, $\mu$-quasiregular on $\C\setminus E$, is actually a constant function. As it is in the $BMO$ case, just $1$-dimensional sets are interesting, because of the Stoilow factorization, together with the fact that $\mu$-quasiconformal mappings with $\mu\in W^{1,2}$ do not distort Hausdorff dimension.\\ 

As we know from Corollary \ref{muzerotozero}, if $E$ is such that $\H^1(E)=0$ then also $\H^1(\phi(E))=0$ whenever $\phi$ is $\mu$-quasiconformal. Thus, also $\gamma(\phi(E))=0$. This shows that zero length sets are removable for bounded $\mu$-quasiregular mappings.\\

Now the following step consists of understanding what happens with sets of positive and finite length. It is well known that those sets can be decomposed as the union of a rectifiable set, a purely unrectifiable set, and a set of zero length (see for instance Mattila \cite[p.205]{M}). Hence, we may study them separately.

\begin{lem}\label{muqcrectifiabledistortion}
Let $\phi:\C\rightarrow\C$ be a planar homeomorphism, such that $\phi,\phi^{-1}\in W^{2,1+\varepsilon}_{loc}(\C)$ for some $\varepsilon>0$. Suppose that $\H^1(E)=0$ if and only if $\H^1(\phi(E))=0$. Then,
$$\Gamma\text{ rectifiable }\hspace{1cm}\Longleftrightarrow\hspace{1cm}\phi(\Gamma)\text{ rectifiable}$$
\end{lem}
\begin{proof}
Since $\Gamma$ is a rectifiable set, then
$$\Gamma\setminus Z=\bigcup_{i=1}^\infty \Gamma_i$$
where $Z$ is a zero length set, and each $\Gamma_i$ is a ${\cal C}^1$ curve with
nonsingular points (i.e. with nonzero tangent vector at each point). Thus, there is no restriction in assuming that $\Gamma$ is a ${\cal C}^1$ regular curve. In other words, from now on we will suppose that $\Gamma=\{\alpha(t); t\in(0,1)\}$ for some ${\cal C}^1$ function
$$\alpha:(0,1)\rightarrow\C$$
such that $\alpha'(t)\neq 0$ for each $t\in(0,1)$.\\
Since $\phi\in W^{2,1+\varepsilon}_{loc}$, then $\phi$ is strongly differentiable $C_{1,1+\varepsilon}$-almost everywhere (see for instance \cite{Dor}), and the same happens to $\phi^{-1}$. Thus, the set
$$B=\{z\in \Gamma:\phi\text{ is differentiable at }z\}$$
is such that $C_{1,1+\varepsilon}(\Gamma\setminus B)=0$. In particular, $\H^1(\Gamma\setminus B)=0$. Moreover, since also $\phi^{-1}\in W^{2,1+\varepsilon}_{loc}$, we can apply the chain rule and from $\phi^{-1}\circ\phi(z)=z$ it follows that
$$D\phi^{-1}(\phi(z))=(D\phi(z))^{-1}$$
for $C_{1,1+\varepsilon}$-almost every $z\in B$. Thus we may assume
that for each $z\in B$ we also have  $J(z,\phi)\neq 0$. Notice that
 $\H^1(\phi(\Gamma)\setminus\phi(B))=0$ and also that $\H^1(\phi(\Gamma))$ is $\sigma$-finite.\\
Fix a point $w_0=\phi(z_0)$, where $z_0\in B$, and put $z_0=\alpha(t_0)$ for
some $t_0\in(0,1)$. Then, $\alpha'(t_0)=v\neq 0$. We will show that
$$V=\langle D\phi(z_0)\cdot\alpha'(t_0)\rangle$$
is a tangent line to $\phi(\Gamma)$ at $w_0$. Since this argument works at $\H^1$-almost every $w_0\in\phi(\Gamma)$, we will obtain that $\phi(\Gamma)$ is rectifiable \cite[p. 214, Remark 15.22]{M}. Thus, what we have to show is that for every number $s\in(0,1)$, there is $r>0$ such that
$$\phi(\Gamma)\cap D(w_0, r)\subset X(w_0,V,s).$$

Here $X(w_0,V,s)$ is the cone of center $w_0$, direction $V$ and amplitude $s$, that is,
$$X(w_0,V, s)=\{w\in\C: d(w-w_0, V)<s\,|w-w_0|\}.$$

First, by the chain rule, the function
$\tilde{\alpha}=\phi\circ\alpha$ is differentiable a $t_0$. Thus,
$$
\lim_{t\to t_0}\frac{|\tilde{\alpha}(t)-\tilde{\alpha}(t_0)-{\tilde{\alpha}'(t_0)}\,(t-t_0)|}{|t-t_0|}=0.
$$
Hence, 
$$
\lim_{t\to t_0}\frac{|\tilde{\alpha}(t)-\tilde{\alpha}(t_0)|}{|t-t_0|}=|\tilde{\alpha}'(t_0)|.
$$
Moreover, one has $\tilde{\alpha}'(t_0)=D\phi(z_0)\cdot\alpha'(t_0)$. On the other hand, since $\alpha$ is a regular curve, for each $r_1>0$ there is
$r_2>0$ such that
$$|t-t_0|<r_1\,\,\,\Leftrightarrow\,\,\,|\alpha(t)-z_0|<r_2.$$
Put $z=\alpha(t)$. Then $|t-t_0|<r_1$ if and only if $z\in \Gamma\cap D(z_0, r_2)$. But if $r_1$ is chosen small enough,
$$\aligned
d(\phi(z)-w_0,V)&=\inf_{\lambda\in\R} |\phi(z)-\phi(z_0)-D\phi(z_0)\cdot\alpha'(t_0)\,\lambda|\\
&\leq|\phi(z)-\phi(z_0)-D\phi(z_0)\cdot\alpha'(t_0)\,(t-t_0)|\\
&=\frac{|\tilde{\alpha}(t)-\tilde{\alpha}(t_0)-{\tilde{\alpha}'(t_0)}\,(t-t_0)|}{|t-t_0|}\,\frac{|t-t_0|}{|\tilde{\alpha}(t)-\tilde{\alpha}(t_0)|}\,|\tilde{\alpha}(t))-\tilde{\alpha}(t_0))|\\
&\leq s\,|\tilde{\alpha}'(t_0)|\,\frac{1}{|\tilde{\alpha}'(t_0)|}\,|\phi(z)-\phi(z_0)|.
\endaligned$$
Hence, for a given $s>0$ there exists $r_0>0$ (just take $r_0=r_2$) such that
$$\phi(\Gamma\cap D(z_0, r_0))\subset X(w_0,V,s).$$
Notice also that $\phi$ is a homeomorphism, so that if $r$ is small enough, then there we can choose $r_0$ such that
$$\phi(\Gamma\setminus D(z_0, r_0))\subset\C\setminus D(w_0,r).$$
Therefore, given $s>0$ we can find two real numbers $r,r_0>0$ for which the set
$$\phi(\Gamma)\cap D(w_0,r)=\phi(\Gamma\cap D(z_0,r_0))\cap D(w_0,r)$$
has all its points in the cone $X(w_0,V,s)$. In other words, given $s>0$ there is $r>0$ such that $\phi(\Gamma)\cap D(w_0,r)\subset X(w_0,V,s).$
\end{proof}

In this lemma, the regularity assumption is necessary. In the following example, due to J. B. Garnett \cite{Ga}, we construct a homeomorphism of the plane that preserves sets of zero length and, at the same time, maps a purely unrectifiable set to a rectifiable set.

\begin{exam}
Denote by $E$ the planar $\frac{1}{4}$-Cantor set. Recall that this set is obtained as a countable intersection of a decreasing family of compact sets $E_N$, each of which is the union of $4^N$ squares of sidelength $\frac{1}{4^N}$, and where every father has exactly $4$ identic children. \\
\begin{figure}[ht]
\begin{center}
\includegraphics{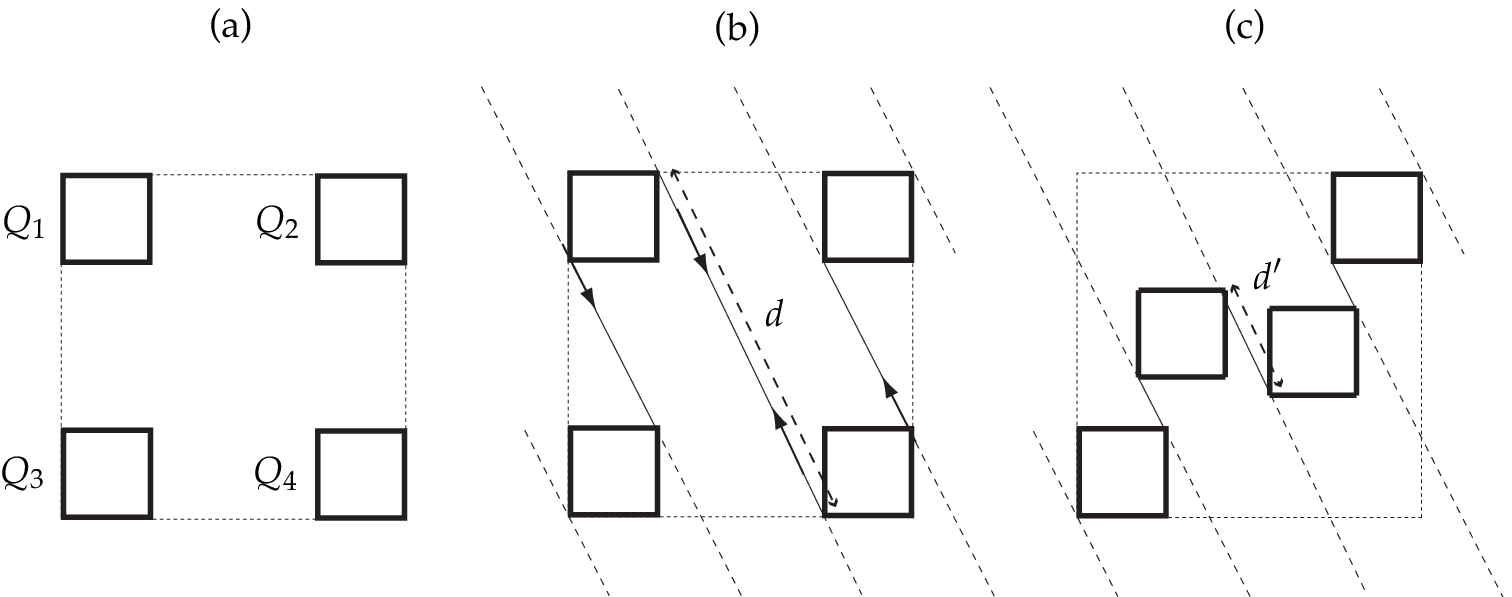}\\
\end{center}
\end{figure}
At the first step, the unit square has $4$ children $Q_1$, $Q_2$, $Q_3$ and $Q_4$. The corners of the squares $Q_j$ are connected with some parallel lines. The mapping $\phi_1$ consists on displacing along these lines the squares $Q_2$ and $Q_3$, while $Q_1$ and $Q_4$ remain fixed. This displacement must be done in such a way that de distance between the images of $Q_2$ and $Q_3$ is positive, since $\phi_1$ must be a homeomorphism.   However, we can do this construction with $d'$ as small as we wish. Our final mapping $\phi$ will be obtained as a uniform limit $\phi=\displaystyle\lim_{N\to\infty}\phi_N$. The other mappings $\phi_N$ are nothing else but copies of $\phi_1$ acting on every one of the different squares in all generations. The only restriction is that the sum of distances $d'$ must be finite. It is clear that this procedure gives a sequence of homeomorphisms $\phi_N$ which converge uniformly to a homeomorphism $\phi$. Further, it can be shown that $\H^1(F)=0$ if and only if $\H^1(\phi(F))=0$. On the other hand, the image of $E$ under the mapping $\phi$ is included in a compact connected set, whose length is precisely the sum of the distances $d'$, which we have chosen to be finite. Therefore, $\phi(E)$ is rectifiable.
\end{exam}

If $\mu\in W^{1,2}$ is a compactly supported Beltrami coefficient, then we know that every $\mu$-quasiconformal mapping belongs to the local Sobolev space $W^{2,q}_{loc}(\C)$ for all $q<2$. Furthermore, we also know that $\phi$ preserves the sets of zero length (even $\sigma$-finite length are preserved), and the same happens to $\phi^{-1}$. Under these hypotheses, we can use for $\phi$ the Lemma \ref{muqcrectifiabledistortion}.

\begin{coro}
Let $\mu\in W^{1,2}$ be a compactly supported Beltrami coefficient, and $\phi$ a $\mu$-quasiconformal mapping.
\begin{enumerate}
\item[(a)] If $E$ is a rectifiable set, then $\phi(E)$ is also rectifiable.
\item[(b)] If $E$ is a purely unrectifiable set, then, $\phi(E)$ is also purely unrectifiable.
\end{enumerate}
\end{coro}
\begin{proof}
The first statement comes from the above lemma. Indeed, in this situation $\mu$-quasiconformal maps preserve sets of zero length and have the needed Sobolev regularity. For the second, let $\Gamma$ be a rectifiable curve. Then,
$$\H^1(\phi(E)\cap \Gamma)=0\Leftrightarrow \H^1(E\cap\phi^{-1}(\Gamma))=0$$
but since $E$ is purely unrectifiable, all rectifiable sets intersect $E$ in a set of zero length. Thus, the result follows.
\end{proof}

\begin{theo}
Let $\mu\in W^{1,2}$ be a compactly supported Beltrami coefficient, and $\phi$ a $\mu$-quasiconformal mapping. Let $E$ be such that $\H^1(E)$ is $\sigma$-finite. Then,
$$\gamma(E)=0\hspace{1cm}\Longleftrightarrow\hspace{1cm}\gamma(\phi(E))=0$$
\end{theo}
\begin{proof}
By Corollaries \ref{muzerotozero} and \ref{mufinitetofinite}, if $\H^1(E)$ is positive and $\sigma$-finite, then $\H^1(\phi(E))$ is positive and $\sigma$-finite. Hence, we may decompose $\phi(E)$ as
$$\phi(E)=\bigcup_n R_n\cup N_n\cup Z_n$$
with $R_n$ rectifiable sets, $N_n$ purely unrectifiable sets, and $Z_n$ zero length sets. Notice that $\gamma(N_n)=0$ because purely unrectifiable sets of finite length are removable for bounded analytic functions \cite{D}, and also $\gamma(Z_n)=0$ since $\H^1(Z_n)=0$. Thus, due to the semiadditivity of analytic capacity \cite{T1}, we get
$$\gamma(\phi(E))\leq C\,\sum_n\gamma(R_n)$$
However, each $R_n$ is a rectifiable set, so that $\phi^{-1}(R_n)$ is also rectifiable. Now, since $E$ has $\sigma$-finite length, the condition $\gamma(E)=0$ forces that $E$ cannot contain any rectifiable subset of positive length, so that $\H^1(\phi^{-1}(R_n))=0$ and hence $\H^1(R_n)=0$. Consequently, $\gamma(\phi(E))=0$.
\end{proof}

The above theorem is an exclusively qualitative result. Therefore, we must not hope for any improvement in a quantitative sense. Namely, in the bilipschitz invariance of analytic capacity by Tolsa \cite{T2}, it is shown that a planar homeomorphism $\phi:\C\to\C$ satisfies $\gamma(\phi(E))\simeq\gamma(E)$ for every compact set $E$ if and only if it is a bilipschitz mapping, while Example \ref{W12example} shows that there exist $\mu$-quasiconformal mappings $\phi$ in the above hypotheses, which are not Lipschitz continuous. 


\section{Beltrami coefficient in $W^{1,p}$, $\frac{2K}{K+1}<p<2$}\label{muW1p}

In this section, we will try to understand the situation when the Beltrami coefficient $\mu$ lies in the Sobolev space $W^{1,p}$, for some $p\in(\frac{2K}{K+1},2)$, where as usually we assume $\|\mu\|_\infty\leq\frac{K-1}{K+1}$. As we showed in Proposition \ref{regularity}, under this assumption every $\mu$-quasiconformal mapping $\phi$ belongs to $W^{2,q}_{loc}$ for each $q<q_0$, where
$$\frac{1}{q_0}=\frac{1}{p}+\frac{K-1}{2K}.$$
Note that we always have $1<q_0<\frac{2K}{2K-1}<2$. We will also denote $p_0=\frac{q_0}{q_0-1}$, so that
$$
\frac{1}{p_0}=\frac{K+1}{2K}-\frac{1}{p}.
$$
Here we always have $p_0\in(2K,\infty)$.\\
\\
Our first goal is to prove an analogous result to Theorem \ref{weylslemmaforbeltrami} (the Weyl's Lemma for the Beltrami operator) in this weaker situation. Let $\Omega$ be either the unit disk $\D$ or the whole plane $\C$. We start by introducing the class of functions
$$E(\Omega)=E^{q,K}(\Omega)=W^{1,q}(\Omega)\cap L^{\frac{2K}{K-1},\infty}(\Omega)$$
for every $1<q<q_0$, equipped with the obvious seminorm
$$\|f\|_E=\|f\|_{\frac{2K}{K-1},\infty}+\|Df\|_q$$
Since ${\cal D}\subset E^{q,K}\subset W^{1,q}$, we get that ${\cal D}$ is dense in $E^{q,K}$.  Further, given $\varphi\in E^{q,K}$,
$$\aligned
\|\mu\varphi\|_E
&=\|\mu\varphi\|_{\frac{2K}{K-1},\infty}+\|D(\mu\varphi)\|_q\\
&\leq\|\mu\|_\infty\,\|\varphi\|_{\frac{2K}{K-1},\infty}+\|D\mu\,\varphi\|_q+\|\mu\|_\infty\,\|D\varphi\|_q\\
&\leq\|\mu\|_\infty\,\left(\|\varphi\|_{\frac{2K}{K-1},\infty}+\|D\varphi\|_q\right)+\|D\mu\|_p\,\|\varphi\|_\frac{pq}{p-q}\\
&\leq\|\mu\|_\infty\,\left(\|\varphi\|_{\frac{2K}{K-1},\infty}+\|D\varphi\|_q\right)+\|D\mu\|_p\,\|\varphi\|_{\frac{2K}{K-1},\infty}
\leq\left(\|\mu\|_\infty+\|D\mu\|_p\right)\,\|\varphi\|_E
\endaligned
$$
In other words, if $1<q<q_0$ then the class $E^{q,K}$ is invariant under multiplication by $\mu$.

\begin{rem}
In the case $p=2$, one has $\|\varphi\|_{\frac{2q}{2-q}} \le \|D\varphi\|_{q}$ for any $q<2$ by the Sobolev embedding. Thus, for any Beltrami coefficient $\mu\in W^{1,2}$ the Sobolev spaces $W^{1,q}$ are invariant under multiplication by $\mu$. However, in our new situation $\frac{2K}{K+1}<p<2$, the Sobolev embedding does not suffice in order to say that for $q<q_0$ the space $W^{1,q}$ is invariant under multiplication by $\mu$.
\end{rem}

The following proposition shows the precise reasons for introducing the class $E^{q,K}$.

\begin{lem}\label{previous}
Let $q\in(1,q_0)$ be fixed, and let $f\in L^\frac{q}{q-1}_{loc}$. Let $\mu\in W^{1,p}$ be a compactly supported Beltrami coefficient. Then, the distribution $\overline\partial f-\mu\partial f$ acts continuously on $E^{q,K}$ functions.
\end{lem}
\begin{proof}
If $f\in L^{\frac{q}{q-1}}$, then clearly $\partial f$ and $\overline\partial f$ act continuously on $W^{1,q}$ functions. However, multiplication by $\mu$ need not be continuous on $W^{1,q}$, so that some extra regularity must be assumed on testing functions. Namely, if $\varphi\in E^{q,K}$ is compactly supported,
$$\aligned
|\langle\overline\partial f-\mu\,\partial f,\varphi\rangle|
&\leq|\langle f,\overline\partial\varphi\rangle|+|\langle f,\partial(\mu\,\varphi)\rangle|\\
&\leq\|f\|_{\frac{q}{q-1}}\,\|\partial\varphi\|_q+\|f\|_{\frac{q}{q-1}}\,\|\partial(\mu\,\varphi)\|_q\\
&\leq\left(1+\|\mu\|_\infty+\|D\mu\|_p\right)\,\|f\|_{\frac{q}{q-1}}\,\|\varphi\|_{E^{q,K}}
\endaligned$$
and the statement follows.
\end{proof}

\begin{theo}
Let $f$ be in $L^{p_0+\varepsilon}_{loc}$ for some $\varepsilon>0$, and assume that
$$\langle\overline\partial f-\mu\,\partial f, \psi\rangle=0$$
for each $\psi\in{\cal D}$. Then, $f$ is $\mu$-quasiregular.
\end{theo}
\begin{proof}
Let $\phi:\C\to\C$ be a $\mu$-quasiconformal mapping, and define $g=f\circ\phi^{-1}$. Clearly $g$ is a locally integrable function. 
Thus we may define $\overline\partial g$ as a distribution and for each $\varphi\in{\cal D}$ we have
$$\aligned
\langle\overline\partial g,\varphi\rangle
&=-\langle g,\overline\partial\varphi\rangle\\
&=-\int g(w)\,\overline\partial\varphi(w)\,dA(w)\\
&=-\int f(z)\,\overline\partial\varphi(\phi(z))\,J\phi(z)\,dA(z)\\
&=-\int f(z)\,\big(\partial\phi(z)\,\overline\partial(\varphi\circ\phi)(z)-\overline\partial\phi(z)\,\partial(\varphi\circ\phi)(z)\big)\,dA(z)
\endaligned$$
This expression makes sense, because $f\in L^{p_0+\varepsilon}$, and both $\phi$ and $\varphi\circ\phi$ belong to $W^{2,q}$ for each $q<q_0$, so that the integrant is an $L^{q}$ function. Now observe the following. On one hand,
$$\aligned
-\int f(z)\,\overline\partial\phi(z)\,\partial(\varphi\circ\phi)(z)dA(z)
&=\langle \partial f,\overline\partial\phi\cdot\,\varphi\circ\phi\rangle + \int f(z)\,\partial\overline\partial\phi(z)\,\varphi\circ\phi(z)\,dA(z)\\
\endaligned$$
and here everything makes sense again, since the function $\partial\phi\cdot(\varphi\circ\phi)$ belongs to the class $E^{q,K}$ for every $q<q_0$, where $\partial f$ acts continuously. Moreover, a similar reasonement gives sense to
$$
-\int f(z)\,\partial\phi(z)\,\overline\partial(\varphi\circ\phi)(z)\,dA(z)=\langle\overline\partial f,\partial\phi\cdot\varphi\circ\phi\rangle + \int f(z)\,\overline\partial\partial\phi(z)\,\varphi\circ\phi(z)\,dA(z)
$$
Thus
$$
\langle\overline\partial g,\varphi\rangle=\langle\overline\partial f,\partial\phi\cdot\varphi\circ\phi\rangle-\langle \partial f,\overline\partial\phi\cdot\,\varphi\circ\phi\rangle.
$$
Now, assume that the Beltrami derivative of $f$ vanishes as a linear functional acting on ${\cal D}$. Then, we get from Lemma \ref{previous} that
$$\langle\overline\partial f-\mu\,\partial f, \psi\rangle=0$$
for every function $\psi$ belonging to $E^{q,K}$, and for any $q<q_0$. 
Since multiplication by $\mu$ is continuous on $E^{q,K}$, the linear functional $\mu\,\partial f$ acts continuously on $E^{q,K}$ functions. Then we can write
$$\langle\overline\partial f,\psi\rangle=\langle\mu\,\partial f,\psi\rangle,$$
or equivalently,
$$\langle\overline\partial f,\psi\rangle=\langle\partial f,\mu\,\psi\rangle,$$
for any $\psi\in E^{q,K}$. In particular, if $\varphi\in{\cal D}$ then the compactly supported function $\psi=\partial\phi\cdot(\varphi\circ\phi)$ belongs to $E^{q,K}$ for every $q<q_0$. Hence,
$$
\langle\overline\partial f,\partial\phi\cdot\varphi\circ\phi\rangle=\langle\mu\,\partial f,\partial\phi\cdot\varphi\circ\phi\rangle=\langle\partial f,\mu\,\partial\phi\cdot\varphi\circ\phi\rangle=\langle\partial f,\overline\partial\phi\cdot\varphi\circ\phi\rangle
$$
and therefore
$$\langle\overline\partial g,\varphi\rangle=\langle\overline\partial f,\partial\phi\cdot\varphi\circ\phi\rangle-\langle \partial f,\overline\partial\phi\cdot\,\varphi\circ\phi\rangle=0$$
and hence $g$ is a holomorphic function. Thus, $f$ is $\mu$-quasiregular.
\end{proof}

Once we know that distributional solutions are strong solutions, also under the weaker assumption $\frac{2K}{K+1}<p<2$, it then follows that some removability theorems can be obtained. The arguments in Section \ref{muqcdistortion} may be repeated to obtain similar estimates for the $BMO$, $VMO$ and $\Lip_\alpha$ problems. In fact, a completely analogous result to Lemma \ref{bmomuremovablesets} holds as well under these weaker assumptions.

\begin{lem}
Let $\frac{2K}{K+1}<p<2$, $E$ be a compact set, and $\mu\in
W^{1,p}$ a Beltrami coefficient, with compact support inside of
$\D$. Suppose that $f$ is $\mu$-quasiregular on $\C\setminus E$,
and $\varphi\in{\cal D}$.
\begin{enumerate}
\item[(a)] If $f\in BMO(\C)$, then
$$|\langle\overline\partial f-\mu\,\partial f,\varphi\rangle|\leq\,C\,\left(1+\|\mu\|_\infty+\|\partial\mu\|_p\right)\,\left(\|\varphi\|_\infty+\|D\varphi\|_\infty\right)\,\|f\|_\ast\,\M^1(E)$$
\item[(b)] If $f\in VMO(\C)$, then
$$|\langle\overline\partial f-\mu\,\partial f,\varphi\rangle|\leq\,C\,\left(1+\|\mu\|_\infty+\|\partial\mu\|_p\right)\,\left(\|\varphi\|_\infty+\|D\varphi\|_\infty\right)\,\|f\|_\ast\,\M^1_\ast(E)$$
\item[(c)] If $f\in\Lip_\alpha(\C)$, then
$$|\langle\overline\partial f-\mu\,\partial f,\varphi\rangle|\leq\,C\,\left(1+\|\mu\|_\infty+\|\partial\mu\|_p\right)\,\left(\|\varphi\|_\infty+\|D\varphi\|_\infty\right)\,\|f\|_\alpha\,\M^{1+\alpha}(E)$$
\end{enumerate}
\end{lem}
\begin{proof}
We repeat the argument in Lemma \ref{bmomuremovablesets}, and consider the function $\delta=\delta(t)$ defined by
$$\delta(t)=\sup_{\diam(D)\leq 2t}\left(\frac{1}{|D|}\int_D|f-f_D|^q\right)^\frac{1}{q}$$
when $0<t<1$, and $\delta(t)=1$ if $t\geq 1$. Here $q=\frac{p}{p-1}$. By construction, for each disk $D\subset\C$ we have
$$\left(\frac{1}{|D|}\int_D|f-f_D|^q\right)^\frac{1}{q}\leq\delta\left(\frac{\diam(D)}{2}\right)$$
Now consider the measure function $h(t)=t\,\delta(t)$. Let $D_j$ be a covering of $E$ by disks, such that
$$\sum_jh(\diam(D_j))\leq\M^h(E)+\varepsilon$$
and consider a partition of unity $\psi_j$ subordinated to the covering $D_j$. Each $\psi_j$ is a ${\cal C}^\infty$ function, compactly supported in $2D_j$, $|D\psi_j(z)|\leq\frac{C}{\diam(2D_j)}$ and $0\leq\sum_j\psi_j\leq 1$ on $\C$. In particular, $\sum_j\psi_j=1$ on $\cup_jD_j$. For every test function $\varphi\in{\cal D}$,
\begin{equation}\label{dossumes2}
-\langle\overline\partial f-\mu\,\partial f,\varphi\rangle=\sum_{j=1}^n\langle f-c_j,\overline\partial(\varphi\,\psi_j)\rangle-\sum_{j=1}^n\langle f-c_j,\partial(\mu\,\varphi\,\psi_j)\rangle.
\end{equation}
An analogous procedure to that in Lemma \ref{bmomuremovablesets} gives
$$\aligned
\left|\sum_{j=1}^n\langle f-c_j,\overline\partial(\varphi\,\psi_j)\rangle\right|\lesssim\left(\M^h(E)+\varepsilon\right)\left(\|\varphi\|_\infty+\|D\varphi\|_\infty\right).\endaligned$$
The other sum in (\ref{dossumes2}) is again divided into two parts,
$$
\left|\sum_{j=1}^n\langle f-c_j,\partial(\mu\,\varphi\,\psi_j)\rangle\right|
\leq\sum_j\int_{2D_j}|f-c_j|\,|\partial\mu|\,|\varphi\,\psi_j|+\sum_j\int_{2D_j}|f-c_j|\,|\mu|\,|\partial(\varphi\,\psi_j)|.$$
The second one, as before,
$$
\sum_j\int_{2D_j}|f-c_j|\,|\mu|\,|\partial(\varphi\,\psi_j)|\lesssim\|\mu\|_\infty\,\left(\M^h(E)+\varepsilon\right)\left(\|\varphi\|_\infty+\|D\varphi\|_\infty\right).
$$
For the first term,
$$\aligned
\sum_j\int|f-c_j|\,|\partial\mu|\,|\varphi\,\psi_j|
&\leq\|\varphi\|_\infty\,\sum_j\left(\int|f-c_j|^q\,|\psi_j|\right)^\frac{1}{q}\,\left(\int|\partial\mu|^p\,|\psi_j|\right)^\frac{1}{p}\\
&\leq\|\varphi\|_\infty\,\left(\sum_j|2D_j|\,\delta(\diam(D_j))^q\right)^\frac{1}{q}\,\left(\sum_j\int|\partial\mu|^p\,\psi_j\right)^\frac{1}{p}\\
&\lesssim\|\varphi\|_\infty\,\left(\sum_j\diam(D_j)^2\,\delta(\diam(D_j))^q\right)^\frac{1}{q}\,\left(\int_\D|\partial\mu|^p\right)^\frac{1}{p}\\
&\lesssim\|\varphi\|_\infty\,\|\partial\mu\|_p\,|\cup_jD_j|^\frac{1}{q}
\endaligned$$
Arguing as in the proof of Lemma \ref{bmomuremovablesets}, this
term turns out to be not the worse one, since area is always
smaller than any other Hausdorff content. The rest of the proof
follows as in Lemma \ref{bmomuremovablesets}.
\end{proof}

Some remarks must be made at this point, just following the ideas of Section \ref{muqcdistortion}. First of all, we observe that the above result has its counterpart in terms of distortion for length. Namely, if $\frac{2K}{K+1}<p<2$ and $\mu$ is a compactly supported Beltrami coefficient, then any $\mu$-quasiconformal mapping $\phi$ satisfies
\begin{equation}\label{zerotozerogeneral}
 \M^1(E)=0\hspace{1cm}\Longrightarrow\hspace{1cm}\M^1(\phi(E))=0
\end{equation}
This was already shown in Section \ref{muqcdistortion} for $\mu\in W^{1,2}$. Now, we are increasing the range of values of $p$ for which this holds and, unexpectedly, this length distortion result is also true for $\frac{2K}{K+1}<p<2$. However, it is not clear a priori if this implication is an equivalence (as it is when $\mu\in W^{1,2}$). Of course, the same can be deduced from the $VMO$ problem. In terms of distortion, the conclusion is that
\begin{equation}\label{finitetofinitegeneral}
 \M^1_\ast(E)=0\hspace{1cm}\Longrightarrow\hspace{1cm}\M^1_\ast(\phi(E))=0
\end{equation}
and again, nothing may be said now on wether this implication is an equivalence. In terms of removability, the analytic and the $\mu$-quasiregular $BMO$ (and $VMO$) problems are the same.\\
\\
For the $\Lip_\alpha$ removability problems, we obtain also the same result that for $\mu\in W^{1,2}$. Although this has no direct implications in terms of distortion of Hausdorff measures, for dimension distortion we get the following.

\begin{coro}
Let $K>1$ and $\frac{2K}{K+1}<p<2$. Let $\mu\in W^{1,p}$ be a compactly supported Beltrami coefficient, and  let $\phi$ be $\mu$-quasiconformal. Then,
$$\dim(E)\leq 1\hspace{1cm}\Longrightarrow\hspace{1cm}\dim(\phi(E))\leq 1$$
\end{coro}
\begin{proof}
If $\dim(E)\leq 1$, then $\H^{1+\alpha}(E)=0$ for all $\alpha>0$, and hence $E$ is $\mu$-removable for $\alpha$-H\"older continuous functions, for every $\alpha>0$. Now let $\beta>0$, and let $h:\C\to\C$ be a $\Lip_\beta$ function, holomorphic on $\C\setminus\phi(E)$. Then, $h\circ\phi$ is a $\Lip_{\beta/K}$ function, $\mu$-quasiregular on $\C\setminus E$ and hence has a $\mu$-quasiregular extension to the whole of $\C$. Then, $h$ extends holomorphically. This means that $\phi(E)$ is removable for holomorphic $\beta$-H\"older continuous functions, so that $\H^{1+\beta}(\phi(E))=0$ \cite{OF}. Since this holds for any $\beta>0$, then we get $\dim(\phi(E))\leq 1$.
\end{proof}

A similar argument proves that
$$\dim(E)\leq t\hspace{1cm}\Rightarrow\hspace{1cm}\dim(\phi(E))\leq 1+K(t-1)$$
for any $t\in (1,1+\frac{1}{K})$.\\
\\
A final remark must be done concerning the question of rectifiability.  As we have said, we do not know if the implications (\ref{zerotozerogeneral}) and (\ref{finitetofinitegeneral}) can be reversed. More precisely, it is not clear if $\mu$-quasiconformal mappings with $\mu\in W^{1,p}$, $\frac{2K}{K+1}<p<2$, preserve sets of zero (or $\sigma$-finite) length. For instance, if $\H^1(E)=0$ then also $\H^1(\phi(E))=0$, but it could happen that $\H^1(E)>0$ and $\H^1(\phi(E))=0$.\\ However, if we restrict ourselves to Beltrami coefficients $\mu\in W^{1,p}$ with $p>\frac{2K^2}{K^2+1}$, then it comes from the Remark after Proposition \ref{muinverse} that the inverse mapping $\phi^{-1}$ is $\nu$-quasiconformal, for some compactly supported Beltrami coefficient $\nu$ satisfying $\|\nu\|_\infty=\|\mu\|_\infty$ and $nu\in W^{1,r}$ for some $r\in(\frac{2K}{K+1}, 2)$. In other words, if $p>frac{2K^2}{K^2+1}$, then $\nu$ belongs to the right range of Sobolev spaces, so that all the comments above apply to $\phi^{-1}$ and we get that both $\phi$ and $\phi^{-1}$ map zero length sets to zero length sets (and the same for $\sigma$-finite length sets). This means that if $\mu\in W^{1,p}$ is a compactly supported Beltrami coefficient, $\frac{2K^2}{K^2+1}<p<2$, and $\phi:\C\to\C$ is $\mu$-quasiconformal, then implications (\ref{zerotozerogeneral}) and (\ref{finitetofinitegeneral}) are actually equivalences. Therefore, $\phi$ is in the situation of Lemma \ref{muqcrectifiabledistortion}, and rectifiable sets are preserved as well under the action of $\phi$. This leads us to the following unexpected result.

\begin{coro}
Let $K>1$, and $\frac{2K^2}{K^2+1}<p<2$. Let $\mu\in W^{1,p}$ be a compactly supported Beltrami coefficient, $\|\mu\|_\infty\leq\frac{K-1}{K+1}$, and let $\phi:\C\to\C$ be a $\mu$-quasiconformal mapping. Then,
$$\gamma(E)=0\hspace*{1cm}\Longleftrightarrow\hspace*{1cm}\gamma(\phi(E))=0$$
for any compact set $E$ with $\sigma$-finite $\H^1(E)$.
\end{coro}

Notice that for small values of $K$, say $K=1+\varepsilon$, it is sufficient to take $p\geq 1+\varepsilon$. That is, given a compactly supported Beltrami coefficient $\mu\in W^{1,1+\varepsilon}$, the conclusion of the above result holds, provided that $\|\mu\|_\infty\lesssim\varepsilon$.

\vskip 1cm
\begin{itemize}
\item[]{Departament de Matem\`atiques, Facultat de Ci\`encies, Universitat Aut\`onoma de Barcelona, 08193-Bellaterra, Barcelona, Catalonia\\
{\it E-mail address:} {albertcp@mat.uab.cat}, {mateu@mat.uab.cat}, {orobitg@mat.uab.cat}}

\item[]{Departamento de Matem\'aticas,
Universidad Aut\'onoma de Madrid,
28049 Madrid,
Spain\\
{\it Email address:} {daniel.faraco@uam.es}}

\item[]{Department of Mathematics and Statistics, University of Jyv\"askyl\"a, P.O.
Box 35 (MaD), FIN-40014, University of Jyv\"askyl\"a, Finland\\
{\it E-mail address:} {zhong@maths.jyu.fi}}
\end{itemize}

\end{document}